\documentclass[preprint,11pt]{elsarticle}
\usepackage[latin1]{inputenc}
\usepackage{amsmath}
\usepackage{amsthm}
\usepackage{amssymb}
\usepackage{amsfonts}
\usepackage{geometry}
\usepackage{graphicx}

\newtheorem{theo}{Theorem}
\newtheorem{cor}[theo]{Corollary}
\newtheorem{prop}[theo]{Proposition}
\newtheorem{pte}[theo]{Properties}
\newtheorem{lemme}[theo]{Lemma}
\newtheorem{conj}[theo]{Conjecture}

\newdefinition{defi}[theo]{Definition}
\newdefinition{rema}[theo]{Remark}
\newdefinition{exe}[theo]{Example}

\newproof{pf}{Proof}

\newcommand{\cA}{\mathcal{A}}
\newcommand{\cM}{\mathcal{M}}
\newcommand{\cS}{\mathcal{S}}
\newcommand{\cR}{\mathcal{R}}
\newcommand{\cU}{\mathcal{U}}
\newcommand{\bcA}{\boldsymbol{\cA}}
\newcommand{\bcS}{\boldsymbol{\cS}}

\newcommand{\VcM}{\Vert \cM_n \Vert}

\newcommand{\bpi}{\boldsymbol{\pi}}
\newcommand{\bmu}{\boldsymbol{\mu}}
\newcommand{\N}{\mathbb{N}}
\newcommand{\Z}{\mathbb{Z}}
\newcommand{\ov}{\overline}
\newcommand{\wh}{\widehat}
\newcommand{\sn}{\sqrt{n}}

\journal{Linear Algebra and its Applications}

\begin{document}

\begin{frontmatter}


\title{Symmetric matrices related to the Mertens function}
\author{Jean-Paul Cardinal}
\address{
Département de Mathématiques\\
Université Paris 13\\
93430 Villetaneuse, France\\
cardinal@math.univ-paris13.fr}

\begin{abstract} 
In this paper we explore a family of congruences over $\N^\ast$ from which one builds a sequence of symmetric matrices related to the Mertens function. From the results of numerical experiments, we formulate a conjecture about the growth of the quadratic norm of these matrices, which implies the Riemann hypothesis. This suggests that matrix analysis methods may come to play a more important role in this classical and difficult problem.
\end{abstract}


\end{frontmatter}

\section{Introduction}

Among the numerous statements equivalent to the Riemann hypothesis, a few have been formulated as matrix problems. The Redheffer matrix $A_n=(a_{i,j})_{1\le i,j \le n}$ is an $n \times n$ matrix defined by $a_{i,j} = 1$ if $j = 1$ or if $i$ divides $j$, and $a_{i,j} = 0$ otherwise. R. Redheffer \cite{Redheffer} proved that the Riemann hypothesis is true if and only if 
$$det(A_n) = O(n^{1/2+\epsilon}), \text{ for every } \epsilon > 0.$$
F. Roesler \cite{Roesler} defined the matrix $B_n=(b_{i,j})_{2\le i,j \le n}$ by $b_{i,j}=i-1$ if $i$ divides $j$, and $b_{i,j}=-1$ otherwise, and he thus proved that the Riemann hypothesis is true if and only if 
$$det(B_n) = O(n!n^{-1/2+\epsilon}), \text{ for every } \epsilon > 0.$$
Both the Redheffer and Roesler matrices are related, via their determinants, to the Mertens function, which is by definition the summatory function of the Möbius function (see \cite{Apostol} p.91). These matrices are not symmetric, and many eigenvalues must be computed to estimate the determinants.\\

The matrix $\cM_n$ that we introduce in this paper is also related to the Mertens function, and it is symmetric.
In Section $2$, we will prove that the Riemann hypothesis is true if 
$$\VcM = O(n^{1/2+\epsilon}), \text{ for every } \epsilon > 0,$$
where $\VcM$ denotes the quadratic norm of the matrix, i.e. its spectral radius. In Section $3$, from the results of a series of numerical experiments, we will conjecture that $\VcM = O(n^{1/2+\epsilon}), \text{ for every } \epsilon > 0$. Consequently, the search for a good estimate of $\VcM$ turns out to be a key problem. At this prospect, it should be noted that, contrary to both the Redheffer and Roesler matrices, only one eigenvalue must be estimated in the matrix $\cM_n$.

\section{Construction of the matrices $\cM_n$}

\subsection{A family of congruences over $\N^\ast$}

\begin{defi}
\label{defR}
To each $n\in \N^\ast$, we associate an equivalence relation $\cR$ over $\N^\ast$, defined by 
$$i \ \cR \ j \Longleftrightarrow \left[ n/i \right]=\left[ n/j \right].$$
\end{defi}

\begin{exe}
For $n=16$, $\cR$ possesses the eight equivalence classes\\
$\lbrace 1\rbrace,\lbrace 2\rbrace,\lbrace 3\rbrace,\lbrace 4 \rbrace,\lbrace 5\rbrace ,
\lbrace 6,7,8\rbrace,\lbrace 9,10,11,12,13,14,15,16 \rbrace$ and $\lbrace 17,18,\cdots\rbrace$.
\end{exe}

Due to the interval structure of these classes, we can unambiguously identify each class by its largest representative (with the convention that $\infty$ denotes the largest representative of the unbounded class).
For clarity, the representatives are written in plain letters, and the classes are written in letters with a hat. We denote by $\ov\cS$ the set of the largest representatives and by $\wh\cS$ the set of the classes, i.e. $\wh\cS=\N^\ast/\cR$.
We also set $\cS=\ov\cS\setminus\lbrace \infty \rbrace$.\\

Throughout Section $2$, most of the objects we define, such as the sets $\cR,\cS$ defined above, and the matrix $\cM$ introduced in Proposition \ref{TUT}, depend on the integer $n$. However, in order to simplify the notations in this section, we will not index these objects by $n$, since there is no risk of ambiguity.

\begin{exe}
\label{exeS}
 For n = 16,\\
$\ov\cS=\{1,2,3,4,5,8,16,\infty \}$, 
$\cS=\{1,2,3,4,5,8,16 \}$,
$\wh\cS=\{\wh{1},\wh{2},\wh{3},\wh{4},\wh{5},\wh{8},\wh{16},\wh{\infty} \}$,\\  
${\wh 1}=\lbrace 1\rbrace,{\wh 2}=\lbrace 2\rbrace,
{\wh 3}=\lbrace 3\rbrace,{\wh 4}=\lbrace 4 \rbrace,{\wh 5}=\lbrace 5\rbrace ,{\wh 8}=\lbrace 6,7,8\rbrace, \wh{16}=\lbrace 9,10,11,12,13,14,15,16 \rbrace$, \\
 $\wh{\infty}=\lbrace 17,18,\cdots \rbrace$.
\end{exe}

\begin{prop}
\label{diezS}
Let $n$ be fixed in $\N^\ast$, and let $\cS$ be defined as above, i.e. $\cS$ is the set of the
 largest representatives of the classes of $\cR$. For each $k$ in $\cS$, we set $\ov{k}=[n/k]$.
\begin{enumerate}
\item For each $k$ in $\cS$, we have $\ov{k}\in\cS$ and $\ov{\ov{k}}=k$, which means that the map 
$k \mapsto \ov{k}$ is a decreasing involution on $\cS$. Actually, $k \mapsto \ov{k}$ is just
the order reversing map on $\cS$.
\item 
The set $\cS$ can be described precisely by the following alternative :
$$
\begin{array}{ll}
\text{If } n<[\sn]^2+[\sn] & \text{then } \cS=\left\lbrace 1,\cdots,[\sn]=\ov{[\sn]},\cdots,n=\ov{1} \right\rbrace ,\\
 & \text{hence } \#\cS=2[\sn]-1; \\
\text{if } n\ge [\sn]^2+[\sn] & \text{then } \cS=\left\lbrace 1,\cdots,[\sn],\ov{[\sn]},\cdots,n=\ov{1} \right\rbrace , \\
& \text{hence } \#\cS=2[\sn].
\end{array}
$$
 \end{enumerate}
\end{prop}

\begin{pf}
\begin{enumerate}

\item
Let $k\in\cS$ and $\ov{k}=[n/k]$, which means that $k\ov{k}\le n < k\ov{k}+k$. That is to say,
$\dfrac{n}{\ov{k}+1} < k \le \dfrac{n}{\ov{k}}$.\\
Since $k$ is an integer, it follows that
$\displaystyle{\left[ \frac{n}{\ov{k}+1}\right] < \left[ \frac{n}{\ov{k}}\right]}$,
which proves that $\ov{k}\in\cS$.\\
Let $k\in\cS$. Since $\left[ \dfrac{n}{k+1}\right] <\left[ \dfrac{n}{k}\right]$, 
it follows that
$ \dfrac{n}{k+1} <\left[ \dfrac{n}{k}\right] =\ov{k}$. That is to say,
$n<k\ov{k}+\ov{k}$. From this last inequality and the fact that $k\ov{k}\le n$, we deduce that 
$\displaystyle{k\le \frac{n}{\ov{k}} <k+1}$, which means that $k=[n/\ov{k}]=\ov{\ov{k}}$.

\item
We begin to prove that each singleton set $\left\lbrace k\right\rbrace $, with
$1\le k <\sn$, is a class. Indeed, if $k<[\sn ]$, then we have successively
 $k+1\le [\sn ]$, $\dfrac{n}{k}-\dfrac{n}{k+1}=\dfrac{n}{k(k+1)}>1$,
and $\left[ \dfrac{n}{k+1}\right] <\left[ \dfrac{n}{k}\right]$.
This last inequality means that $k$ and $k+1$ do not belong to the same class, which is to say, $k\in\cS$.\\
Considering now the case $k=\left[ \sn\right] $, there are two possibilities :
\begin{enumerate}
\item either $\dfrac{n}{k}-\dfrac{n}{k+1}=\dfrac{n}{k(k+1)} \ge 1$, hence
$\left[ \dfrac{n}{k+1}\right] <\left[ \dfrac{n}{k}\right] $,
\item or $\dfrac{n}{k}-\dfrac{n}{k+1}=\dfrac{n}{k(k+1)} < 1$, hence 
$\displaystyle{\frac{n}{k+1}<k\le\frac{n}{k}}$, and since $k$ is an integer,
 it follows that $\left[ \dfrac{n}{k+1}\right] <\left[ \dfrac{n}{k}\right] $.
\end{enumerate}
In both cases, $[\sn]$ and $[\sn]+1$ do not belong to the same class,
so $[\sn]\in\cS$.\\
Now that we have proved that $\left\lbrace 1,2,\cdots,\left[ \sn\right] \right\rbrace \subset \cS$,
let $k\in\cS$ with $k > [\sn]$. 
Therefore, $k$ satisfies the inequalities $k \ge \sn$, $n/k \le \sn$ and $\ov{k} \le [\sn]$.
In other words, $[\sn]$ is the largest element of $\cS$ such that $k\le \ov{k}$. 
Using the fact that $k \mapsto \ov{k}$ is a decreasing involution on $\cS$, and distinguishing the two cases $[\sn] = \ov{[\sn]}$ and $[\sn] < \ov{[\sn]}$, we deduce the expected form of $\cS$.\\
To conclude, we rewrite the condition $[\sn] = \ov{[\sn]}$. Set $\sn=k+\alpha$ with $k\in\N^\ast$ and
$0\le \alpha < 1$.
We have $n = k^2 + 2\alpha k +\alpha^2$, i.e.: $ n/k = k + 2\alpha +\alpha^2/k$, so the following equivalences hold :
$$[\sn] = \ov{[\sn]} \Leftrightarrow 2\alpha +\alpha^2/k < 1 \Leftrightarrow \alpha^2+2k\alpha-k <0
\Leftrightarrow n < k^2 +k ,$$
and this completes the description of $\cS$.

 \end{enumerate}
\end{pf}

\begin{rema}
A synthetic formula for $\#\cS$, one valid for all $n\in\N^\ast$, is 
$$
\#\cS=[\sn]+[\sqrt{n+1/4}-1/2].
$$
\end{rema}

\begin{lemme}
For all positive integers $n,i,j$, we have 
$$
\left[\left[ n/i \right]/j\right]=\left[ n/ij \right].
$$
\end{lemme}

\begin{pf}
 Write $n/i=u+\alpha$, with $u\in\N$ and $0\le \alpha <1$.\\ Hence
$[n/i]/j=u/j$ and $n/ij=u/j+\alpha/j$, from which it follows that
$\left[\left[ n/i \right]/j\right]\le \left[ n/ij \right]$.
If this last inequality were strict, then there would exist an integer $v$ such that
  $u/j < v \le u/j+\alpha/j$.
Hence, $u < vj \le u+\alpha < u+1$, which is impossible since both $u$ and $vj$ are integers.
\end{pf}

\begin{prop}
\label{multiR}
 $\cR$ is compatible with the multiplication over $\N^\ast$, meaning that, for all $i,j,k \in \N^\ast$, we have $i \ \cR \ j \Longrightarrow ik \ \cR \ jk$. Therefore, the formula $\wh i\wh j=\wh{ij}$ defines an induced multiplication over $\wh\cS$ (recall that $\wh i$ denotes the class of $i$, and $\wh \infty$ the class of every integer strictly larger than $n$).
\end{prop} 

\begin{pf}
Assume that $i \ \cR \ j$, so $\left[ n/i \right]=\left[ n/j \right]$. 
Using the previous lemma we deduce \\
$[ n/ik ]=\left[ [n/i]/k \right] =\left[ [n/j]/k \right] =[ n/jk]$,
which is to say, $ik \ \cR \ jk $.
\end{pf}

The set $\N^\ast$, equipped with the usual multiplication, is a commutative semigroup. Since $\cR$ is compatible with the multiplication, the quotient set $\wh\cS=\N^\ast/\cR$, equipped with the induced multiplication, is also a commutative semigroup. From now on, for each fixed $n$ in $\N^\ast$, the equivalence relation $\cR$ will be called a congruence, and any two integers $i,j$ such that $i\cR j$ will be said to be congruent.
 
\begin{exe}
Multiplication table of $\wh\cS=\N^\ast/\cR$, for $n=16$.

\begin{footnotesize}
$$
\begin{array}{c|cccccccc}
    &   \wh{1} &   \wh{2} &   \wh{3} &   \wh{4} &   \wh{5} &   \wh{8} &   \wh{16} &  \wh{\infty} \\ \hline
  \wh{1} &   \wh{1} &   \wh{2} &   \wh{3} &   \wh{4} &   \wh{5} &   \wh{8} &   \wh{16} & \wh{\infty} \\ 
  \wh{2} &   \wh{2} &   \wh{4} &   \wh{8} &   \wh{8} &   \wh{16} &   \wh{16} &  \wh{\infty} &  \wh{\infty} \\ 
  \wh{3} &   \wh{3} &   \wh{8} &   \wh{16} &   \wh{16} &   \wh{16} &  \wh{\infty} &   \wh{\infty} &  \wh{\infty}\\
  \wh{4} &   \wh{4} &   \wh{8} &   \wh{16} &   \wh{16} &  \wh{\infty} &  \wh{\infty} &  \wh{\infty} &  \wh{\infty}\\
  \wh{5} &   \wh{5} &   \wh{16} &   \wh{16} &  \wh{\infty} &  \wh{\infty} &  \wh{\infty} &  \wh{\infty} &  \wh{\infty} \\
  \wh{8} &   \wh{8} &   \wh{16} &  \wh{\infty} &  \wh{\infty} &  \wh{\infty} &  \wh{\infty} &  \wh{\infty} &  \wh{\infty} \\
  \wh{16} &   \wh{16} &  \wh{\infty} &  \wh{\infty} &  \wh{\infty} &  \wh{\infty} &  \wh{\infty} &  \wh{\infty} &  \wh{\infty} \\
\wh{\infty} & \wh{\infty} &  \wh{\infty} &  \wh{\infty} &  \wh{\infty} &  \wh{\infty} &  \wh{\infty} &  \wh{\infty} &  \wh{\infty} 
\end{array} 
$$
\end{footnotesize}
\end{exe}

\subsection{Three $\Z$-algebras}
The $\Z$-algebra of a semigroup $G$ is the set $\Z^G$ of maps from $G$ to $\Z$ equipped with the convolution product $\star$, defined naturally as follows : if $a$ and $b$ are elements of $\Z^G$, 
then $c=a\star b$ is the map defined by
$$\forall t\in G,\ c\left( t\right)=\sum_
{
r,s\in G : rs=t}
a(r)b(s).$$
Of course, this makes sense if the above sum is finite, a condition that is always satisfied in the remainder of this paper. 

\subsubsection{The algebra $A$ of the semigroup $\N^\ast$}
The algebra ${A}=\Z^{\N^\ast}$ of the semigroup $\N^\ast$ is the algebra of Dirichlet series (with integer coefficients) equipped with the convolution product $\star$, also called Dirichlet product (see \cite{Apostol} p.29). This algebra possesses some well-known properties :

\begin{pte}
\label{dirichlet}
\begin{enumerate}
\item
If ${ a}=(a_1,\cdots,a_k,\cdots)$ and ${ b}=(b_1,\cdots,b_k,\cdots)$ are elements of ${A}$,
and ${ c}=a\star b$, then $c_k=\sum_{ij= k} a_i b_j$.
In particular, for $e_i,e_j$, the $i$-th and $j$-th vectors, respectively, of the canonical basis of ${A}$, then we have $e_i\star e_j= e_{ij}$.
\item
The unit of $A$ is $e_1=(1,0,\cdots,0,\cdots)$. An element ${ a}=(a_1,\cdots,a_n,\cdots)$ is invertible if and only if $a_1 = \pm 1$. 
\item
The inverse of $u=(1,1,\cdots,1,\cdots)$ is $\mu$, the Möbius sequence (see \cite{Apostol} p.31).
\end{enumerate}
\end{pte}

\subsubsection{The algebra $\wh A$ of the semigroup $\wh\cS$, the quotient algebra $\bcA$}
If we consider the semigroup $\wh\cS={\N^\ast}/\cR$ equipped with the induced product, then its algebra ${\wh A}=\Z^{\wh\cS}$ is a $\Z$-algebra of dimension $\#\wh\cS$, for which $\wh\cS$ is a basis.
However, more interesting for our purposes is the quotient algebra
$\bcA = \wh A /{\wh\infty{\wh A}}$ where $\wh\infty{\wh A} = \Z\wh\infty$ is the principal ideal of $\wh A$ generated by $\wh\infty$.
Let $\varpi$ be the canonical projection of ${\wh A}$ onto $\bcA$. We print in bold the images under $\varpi$ of the vectors of the basis $\wh\cS$, and more generally any vector in $\bcA$. For instance,
$\varpi(\wh k)={\bf k}$, $\varpi(\wh\infty)={\bf 0}$.
As $k$ runs through the set $\cS$ (see Example \ref{exeS}), $\bf k$ runs through a set denoted by $\bcS$, which is a basis of $\bcA$, called the canonical basis of $\bcA$. Of course $\#\bcS=\#\cS$, a quantity that has been computed in Proposition \ref{diezS}.
Using these notations, it is easy to construct the multiplication table of the basis $\bcS$, from the multiplication table of $\wh \cS$, as follows :
\begin{enumerate}
\item
Remove the last row and the last column from the table of $\wh \cS$, and replace the remaining symbols $\wh\infty$ by ${\bf 0}$ (this expresses the fact that $\varpi(\wh\infty)={\bf 0}$).
\item
Remove the hats and rewrite the integers in bold letters (this expresses the fact that $\varpi(\wh k)={\bf k}$).
\end{enumerate}

\begin{exe}
Multiplication table of the canonical basis $\bcS$, for $n=16$.\\

\begin{footnotesize}
\boldmath
$$
\begin{array}{c|ccccccc}
    &   1 &   2 &   3 &   4 &   5 &   8 &   16  \\ \hline
  1 &   1 &   2 &   3 &   4 &   5 &   8 &   16  \\ 
  2 &   2 &   4 &   8 &   8 &   16 &   16 &  0  \\ 
  3 &   3 &   8 &   16 &   16 &   16 &  0 &   0 \\
  4 &   4 &   8 &   16 &   16 &  0 &  0 &  0 \\
  5 &   5 &   16 &   16 &  0 &  0 &  0 &  0  \\
  8 &   8 &   16 &  0 &  0 &  0 &  0 &  0  \\
  16 &   16 &  0 &  0 &  0 &  0 &  0 &  0  
\end{array} 
$$
\end{footnotesize}
\end{exe}

\subsection{A natural morphism from $A$ to ${\bcA}$}

\begin{defi}
We call $\vartheta$ the morphism of $\Z$-modules defined by 
$$
\vartheta :
\left\vert
\begin{array}{c}
A \mapsto \wh A \\
e_i \mapsto \wh i
\end{array} 
\right. \text{, where } e_i \text{ denotes the $i$-th vector of the canonical basis of } A.
$$
\end{defi}

\begin{prop}
$\vartheta$ is a morphism of $\Z$-algebras. That is, for all $a,b$ in $A$, we have 
$$
\vartheta(a\star b)=\vartheta(a)\vartheta(b).
$$
\end{prop}

\begin{pf}
Since $\vartheta$ is linear, we have only to prove the result for arbitrary basis vectors $e_i,e_j$.
 Using Propositions \ref{multiR} and \ref{dirichlet}, we have
$\vartheta(e_i\star e_j)=\vartheta(e_{ij})=\wh{ij}=\wh i\wh j=\vartheta(e_i)\vartheta(e_j)$.
\end{pf}

\begin{prop}
\label{bpi_morph}
The map $\bpi=\varpi \circ \vartheta$ from $A$ to $\bcA$ is a morphism of $\Z$-algebras.\\
The image of $a=(a_1,\cdots,a_\kappa,\cdots)\in A$ by the morphism $\bpi$ is
$$
\bpi(a) = \sum_{ k\in \cS}\left( \sum_{\kappa\in {\wh k}}a_\kappa\right) {\bf k}.
$$
\end{prop}

\begin{pf}
The only thing to check is that the sum $\sum_{\kappa\in {\wh k}}a_\kappa$ is well defined.
This results from the fact that, for every $k\in \cS$, $\wh k$ is a finite subset of $\N^\ast$.
\end{pf}

\begin{cor}
\label{mertens}
The morphism $\bpi$ has the following properties : 
\begin{enumerate}
\item
It maps $u=(1,1,\cdots,1,\cdots)$ to ${\bf u}=(\#{\wh k})_{k\in \cS}$.
  If for each $k\in \cS$
 we denote by $k^-$ the predecessor of $k$ in $\cS$, with the convention that $1^-=0$, then
  $\#{\wh k}=k-k^-$. Hence, 
$$
{\bf u}=(k-k^-)_{k\in\cS}.
$$
\item
It maps $\mu=(\mu(1),\mu(2),\cdots,\mu(k),\cdots)$ to
${\bmu}=\left(\sum_{\kappa\in {\wh k}}\mu(\kappa)\right) _{ k\in \cS}$. Hence,
$$
\bmu=\left( M(k)-M(k^-)\right) _{k\in\cS},
$$
where $M$ denotes the Mertens function.
\end{enumerate}
\end{cor}

\begin{exe}
\label{mertens_exe}
For $n=16$ we have
\begin{enumerate}
\item
${\bf u}=(1,1,1,1,1,3,8)={\bf 1} + {\bf 2} + {\bf 3} + {\bf 4} + {\bf 5} + 3\ {\bf 8} + 8\ {\bf 16}$.
\item
${\bmu}=(1,-1,-1,0,-1,0,1)={\bf 1-2-3-5+16}\\ 
=\mu(1){\bf 1}+\mu(2) {\bf 2}+\mu(3) {\bf 3}+\mu(4) {\bf 4}+\mu(5) {\bf 5}
+\left( M(8)-M(5)\right)  {\bf 8}+\left( M(16)-M(8)\right) {\bf 16}$.
\end{enumerate}
\end{exe}

\subsection{The regular representation of the algebra $\bcA$}
For every ${\bf a} \in \bcA$, the map
$$
\label{eq:linear_rep}
\left\vert
\begin{array}{c}
\bcA \mapsto \bcA \\ 
{\bf x} \mapsto {\bf a  x}
\end{array} 
\right.
$$
is linear, and it is represented in the canonical basis $\bcS$ of $\bcA$ by a matrix $\rho({\bf a})$. We denote by $s=\#\bcS$ the dimension of $\bcA$ ($s$ was computed as a function of $n$ in Proposition $\ref{diezS}$), and we denote by $\cM_s(\Z)$, the algebra of square matrices of size $s$ with integer entries. 
The map
$$
\label{eq:linear_iso}
\rho:
\left\vert
\begin{array}{c}
{\bcA}\mapsto \mathcal{M}_s(\Z) \\ 
{\bf a} \mapsto \rho({\bf a})
\end{array} 
\right.
$$
is called the regular representation of $\bcA$ (see \cite{James} p.56). Moreover, this representation is faithful, i.e. the morphism $\rho$ is injective.
The set of all matrices $\rho({\bf a})$, for ${\bf a} \in \bcA$, is therefore a commutative sub-algebra  of $\cM_s(\Z)$, of dimension $s$, for which the matrices $\rho({\bf k}),\ {\bf k}\in\bcS$, form a basis. Finally, since there is a natural bijection between $\cS$ and $\bcS$, we choose $\cS$ as the indexing set for the rows and the columns of the matrices $\rho({\bf a})$. For instance, when $n=16$, the last column of any matrix $\rho({\bf a})$ does not have index $7$, but index $16$.

\begin{exe}
\label{representants}
For n=16, in addition to $\rho({\bf 1})$, which is the identity matrix, the matrices representing ${\bf 2, 3, 4, 5, 8, 16}$
(where most of the zero entries are left blank for legibility) are
\begin{center}
\begin{footnotesize}
$
\begin{array}{c|ccccccc|}
\rho({\bf 2}) & { 1} & { 2} & { 3} & { 4} & { 5} & { 8} & { {16}}\\ \hline
{ 1} & & & & & & &\\ 
{ 2} & 1 & & & & & &\\ 
{ 3} & & & & & & &\\ 
{ 4} & & 1 & & & & &\\ 
{ 5} & & & & & & &\\ 
{ 8} & & & 1 & 1 & & &\\ 
{ {16}} & & & & & 1 & 1 &\\ \hline
\end{array} 
$
$
\begin{array}{c|ccccccc|}
\rho({\bf 3}) & { 1} & { 2} & { 3} & { 4} & { 5} & { 8} & { {16}}\\ \hline
{ 1} & & & & & & &\\ 
{ 2} & & & & & & &\\ 
{ 3} & 1 & & & & & &\\ 
{ 4} & & & & & & &\\ 
{ 5} & & & & & & &\\ 
{ 8} & & 1 & & & & &\\ 
{ {16}} & & & 1 & 1 & 1 & &\\ \hline
\end{array} 
$
$
\begin{array}{c|ccccccc|}
\rho({\bf 4}) & { 1} & { 2} & { 3} & { 4} & { 5} & { 8} & { {16}}\\ \hline
{ 1} & & & & & & &\\ 
{ 2} & & & & & & &\\ 
{ 3} & & & & & & &\\ 
{ 4} & 1 & & & & & &\\ 
{ 5} & & & & & & &\\ 
{ 8} & & 1 & & & & &\\ 
{ {16}} & & & 1 & 1 & & &\\ \hline
\end{array} 
$\\
$
\begin{array}{c|ccccccc|}
\rho({\bf 5}) & { 1} & { 2} & { 3} & { 4} & { 5} & { 8} & { {16}}\\ \hline
{ 1} & & & & & & &\\ 
{ 2} & & & & & & &\\ 
{ 3} & & & & & & &\\ 
{ 4} & & & & & & &\\ 
{ 5} & 1 & & & & & &\\ 
{ 8} & & & & & & &\\ 
{ {16}} & & 1 & 1 & & & &\\ \hline
\end{array} 
$
$
\begin{array}{c|ccccccc|}
\rho({\bf 8}) & { 1} & { 2} & { 3} & { 4} & { 5} & { 8} & { {16}}\\ \hline
{ 1} & & & & & & &\\ 
{ 2} & & & & & & &\\ 
{ 3} & & & & & & &\\ 
{ 4} & & & & & & &\\ 
{ 5} & & & & & & &\\ 
{ 8} & 1& & & & & &\\ 
{ {16}} & & 1& & & & &\\ \hline
\end{array} 
$
$
\begin{array}{c|ccccccc|}
\rho({\bf {16}}) & { 1} & { 2} & { 3} & { 4} & { 5} & { 8} & { {16}}\\ \hline
{ 1} & & & & & & &\\ 
{ 2} & & & & & & &\\ 
{ 3} & & & & & & &\\ 
{ 4} & & & & & & &\\ 
{ 5} & & & & & & &\\ 
{ 8} & & & & & & &\\ 
{ {16}} & 1& & & & & &\\ \hline
\end{array} \ .
$
\end{footnotesize}
\end{center}
Furthermore, following Example \ref{mertens_exe}, the matrices representing ${\bf u}$ and ${\bmu}$ are
\begin{center}
\begin{footnotesize}
$
\begin{array}{c|ccccccc|}
\rho({\bf u}) & { 1} & { 2} & { 3} & { 4} & { 5} & { 8} & { {16}}\\ \hline
{ 1} & 1& & & & & &\\ 
{ 2} & 1& 1& & & & &\\ 
{ 3} & 1& 0& 1& & & &\\ 
{ 4} & 1& 1& 0& 1& & &\\ 
{ 5} & 1& 0& 0& 0& 1& &\\ 
{ 8} & 3& 2& 1& 1& 0& 1&\\ 
{ {16}} & 8& 4& 3& 2& 2& 1& 1\\ \hline
\end{array} 
$
\hspace{1cm}
$
\begin{array}{c|ccccccc|}
\rho({\bmu}) & { 1} & { 2} & { 3} & { 4} & { 5} & { 8} & { {16}}\\ \hline
{ 1} & 1  &   &    &   &   &   &\\ 
{ 2} & -1 & 1 &    &   &   &   &\\ 
{ 3} & -1 &  0 &  1 &   &   &   &\\ 
{ 4} &  0  &-1 &  0  & 1 &   &   &\\ 
{ 5} & -1 &  0 &  0  &  0 & 1 &   &\\ 
{ 8} &  0  & -1& -1 & -1& 0  & 1 &\\ 
{ {16}}& 1  & -1& -2 & -1& -2& -1& 1\\ \hline
\end{array} \ ,
$
\end{footnotesize}
\end{center}
on which we verify that for ${\bf a}\in {\bcA}$, the coefficients of ${\bf a}$ in the basis $\bcS$ appear in the first column of $\rho({\bf a})$, cf. Example \ref{mertens_exe}.
\end{exe}

\begin{defi}
For $n\in\N^\ast$, let $T$ be the symmetric matrix of size $s=\#\cS$ whose entries are all $1$'s above the second diagonal and $0$'s strictly below.
\end{defi}

\begin{exe} For $n=16$, the matrix $T$ is
$$T=
\left[
\begin{array}{ccccccc}
  1& 1& 1& 1& 1& 1 &1\\ 
  1& 1& 1& 1& 1& 1 &\\ 
  1& 1& 1& 1& 1&   &\\ 
  1& 1& 1& 1&  &   &\\ 
  1& 1& 1&  &  &   &\\ 
  1& 1&  &  &  &   &\\ 
  1&  &  &  &  &   &\\
\end{array} 
\right].
$$
\end{exe}

\begin{lemme}
\label{TM_sym} 
For every ${\bf k}\in \bcS$ and every $i,j\in\cS$, the following equivalence holds : 
$$\left( T\rho({\bf k}) \right)_{i,j}=1 \Leftrightarrow  ij\le [n/k]\Leftrightarrow  k\le [n/ij].$$
\end{lemme}

\begin{pf}
Let ${\bf k}\in \bcS,\ j\in\cS$, and let $v$ be the column of index $j$ of the matrix $\rho(\bf k)$.
The only non-zero entry of $v$, which is $1$, is located at the index ${l}$ such that ${\bf l = kj}$,
i.e. $[n/l]=[n/jk]$. Therefore, the column $Tv$ is the column of index $l$ of $T$,
which is composed of $1$'s for all indices $i$ such that $i\le \ov{l}=[n/l]$,
and of $0$'s below. Moreover, the integers $l$ and $\ov{l}$ are in symmetric positions in the list $\cS$ (see Proposition \ref{diezS}). From this we deduce that 
$$
 \begin{array}{cccccccc}
&\left( T\rho({\bf k}) \right)_{i,j}=1 & \Leftrightarrow & (Tv)_{i}=1 & \Leftrightarrow & i\le [n/l] 
& \Leftrightarrow & i\le [n/jk] \\
\Leftrightarrow & i\le n/jk & \Leftrightarrow & ij\le n/ k & \Leftrightarrow & ij\le [n/ k] & \Leftrightarrow & k\le [n/ ij] 
\end{array}.$$
\end{pf}

\begin{prop}
For all $\bf{a}\in\bcA$, the matrix $T\rho({\bf a})$ is symmetric.
\end{prop}

\begin{pf}
From Lemma \ref{TM_sym}, the matrices $T\rho({\bf k})$ are symmetric. Moreover, these matrices form a basis of $\bcA$. Therefore, by linearity, the matrix $T\rho({\bf a})$ is symmetric for every $\bf{a}\in\bcA$.
\end{pf}

\begin{prop}
\label{TUT}
If we introduce the notations
$\cU=T\rho({\bf u})$ and $\cM=T\rho({\bmu})$, then the matrices $\cU$ and $\cM$ are symmetric and satisfy the relation
$$
\cM=T\cU^{-1}T.
$$
\end{prop}

\begin{pf}
From Proposition \ref{bpi_morph} and item $\mathit{3}$ of Proposition \ref{dirichlet}, we have $\bmu={\bf u}^{-1}$, from which it follows that
$$\cM=T\rho(\bmu)=T\rho({\bf u})^{-1}=T\left( \cU^{-1}T\right). $$
\end{pf}

\subsection{The matrix $\cM$ consists of values of the Mertens function}

Proposition \ref{TUT} establishes a relation between the matrices $\cU$ and $\cM$. In essence, we can compute $\cM$ by inverting $\cU$ and multiplying the result on both sides by $T$. The following is another relation between $\cU$ and $\cM$, involving the Mertens function.

\begin{prop}
\label{TuTmu}
For every $n\in\N^\ast$, the matrices $\cU$ and $\cM$ can be computed by the formulas
$$
\cU =\left( [n/ij]\right) _{i,j\in \cS}  \text{ and }
\cM =\left( M\left( [n/ij]\right) \right) _{i,j\in \cS},
$$
where $M$ denotes the Mertens function. In other words, $\cM$ can be computed by applying the Mertens function term by term to the matrix $\cU$.
\end{prop}

\begin{pf}
We noted in Corrolary \ref{mertens} that ${\bf u}=\sum_{k\in \mathit{S}}(k-k^-){\bf k}$.
Therefore, by linearity,
$$ \begin{array}{ccc}
 \cU & = & \sum_{k\in \mathit{S}}(k-k^-)T\rho({\bf k}),\\
\cU_{i,j} & = & \sum_{k\in \mathit{S}}(k-k^-)\left( T\rho({\bf k}) \right)_{i,j},
 \end{array}$$
and from Lemma \ref{TM_sym}, we deduce that
$$\cU_{i,j}=\sum_{k\in \mathit{S},\ k\le [n/ij]}(k-k^-)=[n/ij].$$
This completes the proof concerning  $\bf u$.\\
Similarly, from Corollary \ref{mertens}, we have
$\bmu=\sum_{k\in \mathit{S}}\left( M(k)-M(k^-)\right) {\bf k}$.
Hence,
$$
 \begin{array}{ccc}
\cM & = & \sum_{k\in \mathit{S}}\left( M(k)-M(k^-)\right) T\rho({\bf k}),\\
\cM_{i,j} & = & \sum_{k\in \mathit{S}}\left( M(k)-M(k^-)\right) \left( T\rho({\bf k}) \right)_{i,j},\\
\cM_{i,j} & = & \sum_{k\in \mathit{S},\ k\le [n/ij]}\left( M(k)-M(k^-)\right), \\
\cM_{i,j} & = & M\left( [n/ij]\right).
\end{array}
$$
\end{pf}

\begin{exe}
\label{TuTmu_exe}
For $n=16$, the matrices $\cU$ and $\cM$ are
\begin{center}
\begin{footnotesize}
$
\begin{array}{c|ccccccc|}
\cU & { 1} & { 2} & { 3} & { 4} & { 5} & { 8} & { {16}}\\ \hline
 1 & 16& 8& 5& 4& 3& 2&1\\ 
 2 & 8& 4& 2& 2& 1& 1&\\ 
 3 & 5& 2& 1& 1& 1& &\\ 
 4 & 4& 2& 1& 1& & &\\ 
 5 & 3& 1& 1& & & &\\ 
 8 &  2& 1& & & & &\\ 
16 &  1& & & & & &\\ \hline
\end{array} 
$
\hspace{1cm}
$
\begin{array}{c|ccccccc|}
\cM & { 1} & { 2} & { 3} & { 4} & { 5} & { 8} & { {16}}\\ \hline
 1 & -1& -2& -2& -1& -1& 0&1\\ 
 2 & -2& -1& 0& 0& 1& 1&\\ 
 3 & -2& 0& 1& 1& 1& &\\ 
 4 & -1& 0& 1& 1& & &\\ 
 5 & -1& 1& 1& & & &\\ 
 8 & 0& 1& & & & &\\ 
16 & 1& & & & & &\\ \hline
\end{array} \ .
$
\end{footnotesize}
\end{center}
\end{exe}

Throughout this section, the matrices $\cM$, $T$ and $\cU$ were not indexed by the integer $n$, although these matrices depended on $n$. Henceforth, we will use the notations $T_n$, $\cU_n$, $\cM_n$ instead of $T$, $\cU$, $\cM$, to express the dependence of these matrices on $n$.

\begin{theo}
\label{cond_fondamentale}
The Riemann hypothesis is true if 
$$
\VcM = O(n^{1/2+\epsilon}), \text{ for every } \epsilon > 0.
$$
\end{theo}
\begin{pf}
On the one hand, we derive from Proposition \ref{TuTmu} that $M(n)$ is the $(1,1)$-entry of $\cM_n$, and  we know that for every square matrix $A$, one has $\max \vert A_{i,j} \vert \le \Vert A \Vert$ (see \cite{Golub} p.57). Therefore, we have $\vert M(n) \vert \le \VcM$. On the other hand, Littlewood \cite{Littlewood} proved that the Riemann hypothesis is equivalent to the estimate 
$\ M(n)=O(n^{1/2+\epsilon}), \text{ for every } \epsilon > 0$.
\end{pf}

\begin{rema}
Our hope is that a good estimate of the spectral radius of $\cM_n$, i.e. its largest eigenvalue, could eventually be obtained by investigating the smallest eigenvalue of the inverse matrix $T_n^{-1}\cU_n T_n^{-1}$, whose construction is quite simple. 
\end{rema}

In the next section, we will look at the quantity $\VcM$ experimentally, as $n$ varies through a certain range of integers.

\section{Numerical experiments}
This section presents the results of some experimental computations
concerning the growth of the sequence $\VcM$ as $n$ tends to infinity.

\subsection{Regularity of the sequence $\VcM$}

Figures $1$ and $2$ display the sequences $M(n)/\sn$ and  $\VcM/\sn$ respectively, for $n$ running from $10^3$ to $10^6$, with a step of $10^3$. Figure $2$bis shows the same data as in Figure $2$, but displayed in a window of smaller height.

\begin{figure}[ht]
	\centering
		\includegraphics{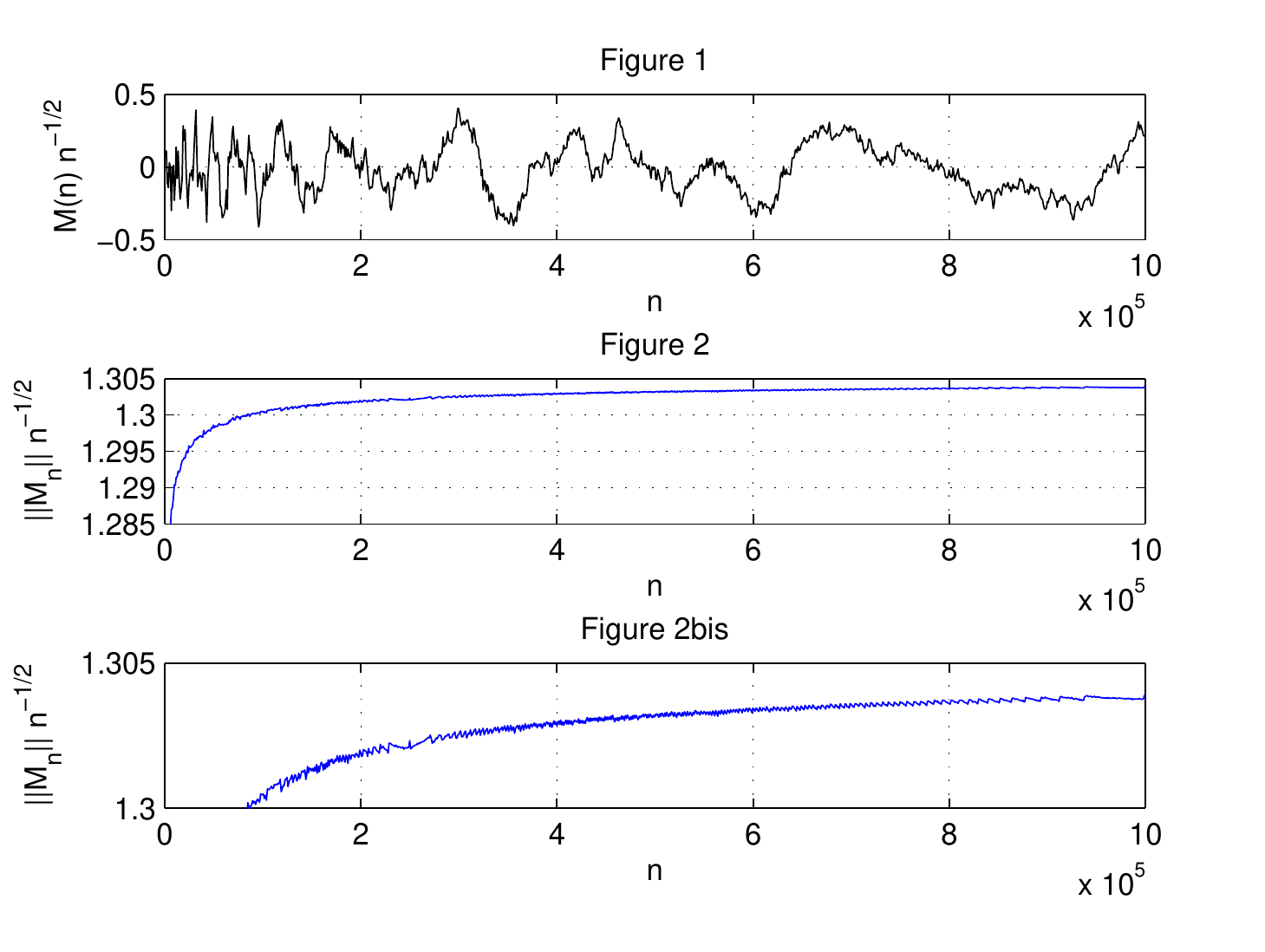}
	\label{fig:figures12}
\end{figure}

We observe that the growth of $\VcM/\sn$ is quite regular, in contrast to the chaotic behavior of $M(n)/\sn$. Not only is the behavior of $\VcM/\sn$ more regular, but Figures $2$ and $2$bis also show that the range in which the sequence $\VcM/\sn$ takes its values is much narrower as $n$ increases, as compared to the case of $M(n)/\sn$. Another important observation is that the growth of $\VcM/\sn$ seems to be relatively slow. We now look more closely at this growth.

\subsection{Experimental convergence of the sequence $\dfrac{\log(\VcM)}{ \log n}$ towards $1/2$}

With regard to Theorem \ref{cond_fondamentale} we now turn our attention to the sequence $w_n=\dfrac{\log(\VcM)}{ \log n}-1/2$.
Figures $3$ to $6$ display the sequence $w_n$ for $n$ taking all the integer values in four intervals centered on the values $n_2=200^2,\ n_3=300^2,\ n_4=400^2$, and $n_5=500^2$.

\begin{figure}[ht]
	\centering
		\includegraphics{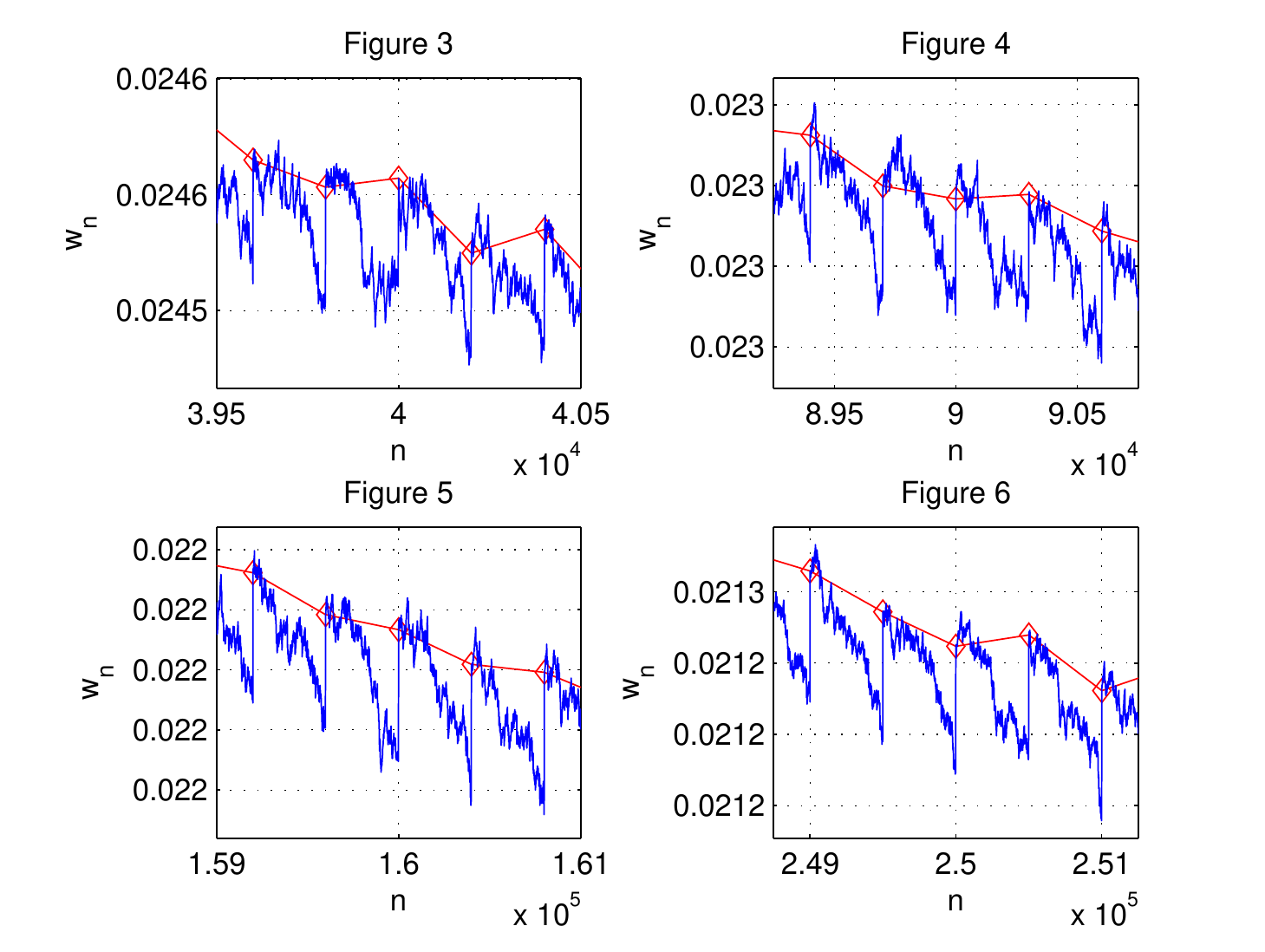}
	\label{fig:figures3456}
\end{figure}

 In these figures, the points $(n,w_n)$ are linked by a blue line, but when $n$ is of the form $k^2$ or $k^2+k$, we plot $(n,w_n)$ as a red diamond. The diamonds are linked by a solid line for greater legibility. We distinguish the two cases because the size of the matrix $\cM_n$ increases by one precisely when $n$ increases from $k^2-1$ to $k^2$ or from $k^2+k-1$ to $k^2+k$. On the figures, this results in small upward jumps in the values of $w_n$. Some structure can be observed in the variations of $w_n$ between every two successive integers $n$ of the form $k^2$ or $k^2+k$. Inside such an interval, the values of $w_n$ seems to follow a random walk of moderate amplitude with a dominant decreasing trend, followed by jumping when $n$ reaches $k^2$ or $k^2+k$. These observations suggest that the overall behavior of the sequence $w_n$ is best described when the values of $n$ are restricted to the forms $k^2$ or $k^2+k$.\\

Figure $7$ displays the sequence $w_n$ for $n$ running from $10^3$ to $10^6$, with the values of $n$ restricted to the forms $k^2$ or $k^2+k$; Figure $8$ shows the same data displayed in loglog axes, i.e. $\log(w_n)$ plotted against $\log(n)$.
\pagebreak

\begin{figure}[ht]
	\centering
		\includegraphics{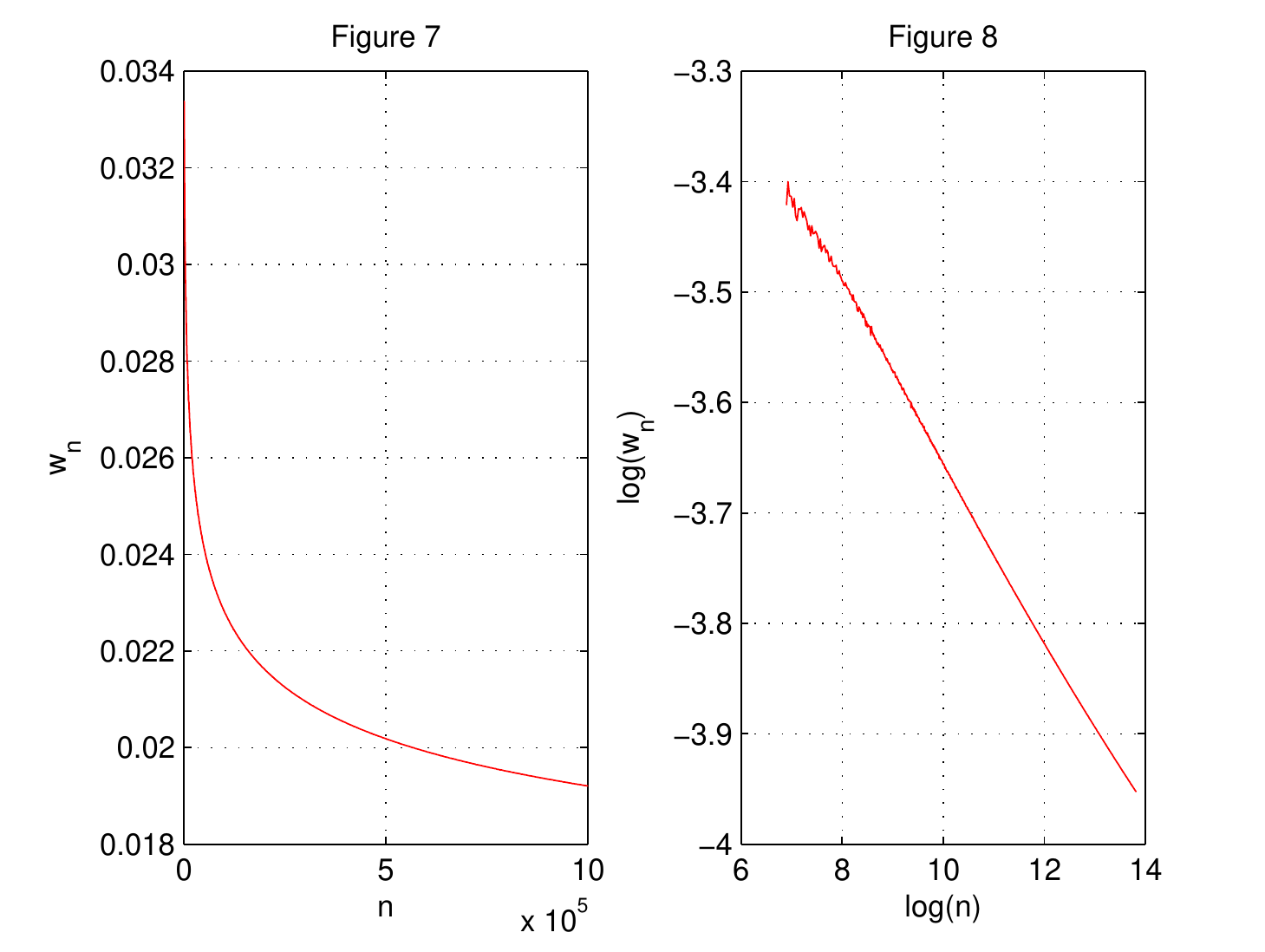}
	\label{fig:figures78}
\end{figure}

We observe in Figure $7$ that $w_n$ is roughly decreasing and remains positive (the positivity results from the fact that $\VcM/\sn\ge 1$ within the range considered), so we may expect that $w_n$ converges to some limit, possibly $0$. The possible convergence of $w_n$ towards $0$ is not very apparent in Figure $7$, but it is more striking in Figure $8$. If, as the graph suggests, this trend were to be confirmed as $n$ increases indefinitely, then the following conjecture would be true :

\begin{conj} 
\label{conj_fondamentale}
We have the estimate :
$$
\VcM =O(n^{1/2+\epsilon}) \text{, for every } \epsilon > 0.
$$
\end{conj}

\section{Conclusion}

We have built a sequence of symmetric matrices $\cM_n$ satisfying $ \vert M(n) \vert \le \VcM$ for all positive integers $n$, where $M$ denotes the Mertens function. Based on numerical evidence, we have conjectured that
$\VcM =O(n^{1/2+\epsilon}), \text{ for every } \epsilon > 0$, a statement which implies the Riemann hypothesis. It may be noted that in no part of this study have we made use of complex variable methods. Finally, the symmetry of the matrices $\cM_n$ suggests that spectral methods in matrix analysis could play a more significant role in the search for a solution to the Riemann hypothesis.

\section{Acknowledgements}
The author wishes to thank Professor Marius Overholt of the University of Troms$\phi$ for his careful reading of this paper and for his most valuable advice.


\begin{thebibliography}{00}

\bibitem{Apostol}
{T.M. Apostol}, {Introduction to Analytic Number Theory}, {Springer}, {1976}.

\bibitem{BFP}
{W.W. Barrett, R.W. Forcade, A.D. Pollington},
{{On the spectral radius of a (0,1) matrix related to the Mertens function}},
{Linear Algebra and its Applications}.
107 (1988)  {151-159}.

\bibitem{BR}
         {W.W. Barrett, D.W. Robinson},
         {{The Jordan 1-structure of a matrix of Redheffer}},
         {Linear Algebra and its Applications}.
         112 (1989) 57-73.
         
         

\bibitem{BJ}
         {W.W. Barrett, T.J. Jarvis},
         {{Spectral properties of a matrix of Redheffer}},
         {Linear Algebra and its Applications}.
         162-64
         (1992)
         673-683.

		
\bibitem{Golub}
		 {G.H. Golub, C.F. van Loan},
		 {Matrix Computations}, {The Johns Hopkins University Press}, {1989}.



\bibitem{Humphries}
         {S.P. Humphries},
         {{Cogrowth of groups and a matrix of Redheffer}},
         {Linear Algebra and its Applications}.
         {265} (1977)
         {101-117}.


\bibitem{James}
		 {G. James, M. Liebeck},
		 {Representations and Characters of Groups}, {Cambridge University Press}, {2006}.


\bibitem{Jarvis}
         {T. J. Jarvis},
         {{A Dominant Negative Eigenvalue of a Matrix of Redheffer}},
         {Linear Algebra and its Applications}.
         {142} (1990)  {141-152}.

\bibitem{Kotnik}
		 {T. Kotnik, J. van de Lune},
		 {{Further Systematic Computations on the Summatory Function of the Möbius Function}},
		 {CWI Report MASR0313},
		 (2003) {1-9}.
		 
\bibitem{Kotnik2}
         {T. Kotnik, J. van de Lune},
         {{On the Order of the Mertens Function}},
         {Experimental Mathematics}.
         {13-4}
         (2004)
         {473-481}.


\bibitem{Littlewood}
         {J.E. Littlewood},
         {Sur la distribution des nombres premiers},
         {C. R. Acad. Sci.}.
         {158}
         (1917)
         {1869-1872}.


\bibitem{Mertens}
         {F. Mertens},
         {{Über eine zahlentheoretische Funktion}},
         {Sitzungsber. Akad. Wiss. Wien}.
         {IIa-106}
         (1897)
         {761-830}.


\bibitem{MVR}
	 {H.L. Montgomery, R.C. Vaughan},
	 {Multiplicative Number Theory I: Classical Theory}, {Cambridge University Press}, {2006}.

\bibitem{Odlyzko}
A.M. Odlyzko, H.J.J. te Riele, Disproof of the Mertens Conjecture,
Journal für die reine und angewandte Mathematik.
357 (1985) 138-160.


\bibitem{Pintz}
         {J. Pintz},
         {{An Effective Disproof of the Mertens Conjecture}},
         {Astérisque}.
         {147-148},
         (1987)
         {325-333}.



\bibitem{Redheffer}
R.M. Redheffer,
{{Eine explizit lösbare Optimierungsaufgabe}},
{Internat. Schriftenreihe Numer. Math}.
36 (1977).

\bibitem{Roesler}
         {F. Roesler},
         {{Riemann's Hypothesis as an Eigenvalue Problem}},
         {Linear Algebra and its Applications}.
         {81}
         (1986)
         {153-198}.


\bibitem{RoeslerII}
         {F. Roesler},
         {{Riemann's Hypothesis as an Eigenvalue Problem II}},
         {Linear Algebra and its Applications}.
         {92}
         (1987)
         {45-73}.
         

\bibitem{RoeslerIII}
         {F. Roesler},
         {{Riemann's Hypothesis as an Eigenvalue Problem III}},
         {Linear Algebra and its Applications}.
         {141} (1990)   {1-46}.

         
         
\bibitem{teRiele}
         {H.J.J. te Riele},
         {{The Mertens Conjecture Revisited}},
        {7th Algorithmic Number Theory Symposium}, {Technische Universität Berlin}, {2006}.


\bibitem{VaughanI}
	 {R.C. Vaughan},
	 {{On the eigenvalues of a Redheffer's matrix I}}, in: Proc. Conf. BYU, Marcel Dekker, 1993, pp. 283-296.


\bibitem{VaughanII}
	 {R.C. Vaughan},
	 {{On the eigenvalues of a Redheffer's matrix II}},
	 {J. Austral. Math. Soc. (Series A)}.
	 {60}
	 (1996)
	 {260-273}.


\bibitem{Wilf}
         {H.S. Wilf},
         {{The Redheffer Matrix of a Partially Ordered Set}},
         {Journal of Combinatorics}.
         {11-4} (2004) R10.
  
\end{thebibliography}
\end{document}